# Stochastic resonance and bifurcation of order parameter in a coupled system of underdamped Duffing oscillators

Ruonan Liu, Yanmei Kang*, Yuxuan Fu
School of Mathematics and Statistics, Xi'an Jiaotong University, Xi'an, Shaanxi, 710049, China
Guanrong Chen
Department of Electronic Engineering, City University of Hong Kong, Hong Kong SAR, China
*Corresponding author email: ymkang@xjtu.edu.cn

The long-term mean-field dynamics of coupled underdamped Duffing oscillators driven by an external periodic signal with Gaussian noise is investigated. A Boltzmann-type *H*-theorem is proved for the associated nonlinear Fokker-Planck equation to ensure that the system can always be relaxed to one of the stationary states as time is long enough. Based on a general framework of the linear response theory, the linear dynamical susceptibility of the system order parameter is explicitly deduced. With the spectral amplification factor as a quantifying index, calculation by the method of moments discloses that both mono-peak and double-peak resonance might appear, and that noise can greatly signify the peak of the resonance curve of the coupled underdamped system as compared with a single-element bistable system. Then, with the input signals taken from laboratory experiments, further observations show that the mean-field coupled stochastic resonance system can amplify the periodic input signal. Also, it reveals that for some driving frequencies, the optimal stochastic resonance parameter and the critical bifurcation parameter have a close relationship. Moreover, it is found that the damping coefficient can also give rise to nontrivial non-monotonic behaviors of the resonance curve, and the resultant resonant peak attains its maximal height if the noise intensity or the coupling strength takes the critical value. The new findings reveal the role of the order parameter in a coupled system of chaotic oscillators.

## 1. Introduction

Stochastic resonance (SR) was originally proposed for explaining the periodic recurrence of the warm and the cold climates [Benzi *et al.*, 1981]. Although SR could not be verified in paleoclimatology, it was subsequently verified by various experiments using, for example, Schmidt trigger circuits [Fauve & Heslot, 1983], bistable ring lasers [Mcnamara *et al.*, 1988], nanomechanical systems [Douglass *et al.*, 1993] and mammalian neuronal networks [Badzey & Mohanty, 2005]. Such experiments were proven successful in explaining the SR phenomenon. A special characteristic of SR lies in that a suitable amount of noise can lead to a distinctive enhancement of a weak input component to a nonlinear system [Jung, 1993; Gammaitoni *et al.*, 1998; Cherubini *et al.*, 2017], similarly to the periodic recurrence of the warm and the cold climates. To date, SR has been widely applied to various engineering fields, ranging from signal processing and detection [Lee *et al.*, 2003; Sun & Kwong, 2007; Fu *et al.*, 2018], early fault diagnosis [Vania & Pennacchi, 2004; Leng *et al.*, 2006; Qiao *et al.*, 2017; Ma *et al.*, 2018], energy harvesting [Harne & Wang, 2013], to image processing [Singh *et al.*, 2017]. Fig. 1 shows an application example of SR in image detection.

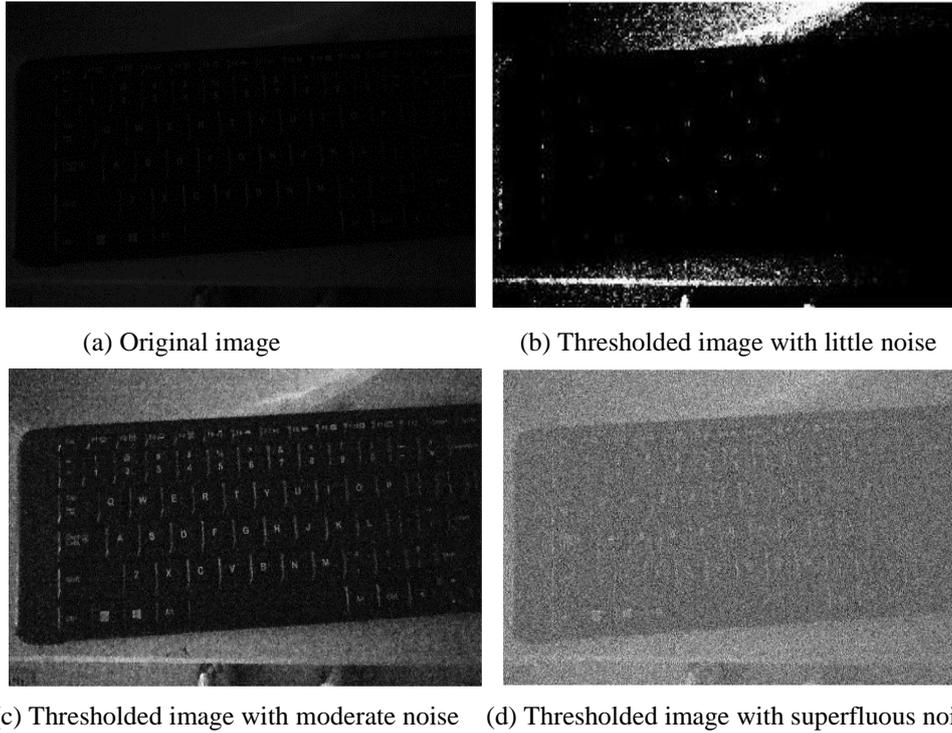

(a) Original image　　　　　　　　(b) Thresholded image with little noise

(c) Thresholded image with moderate noise　　(d) Thresholded image with superfluous noise

Fig.1. Application example of SR in image detection. The original image in (a) is taken in darkness. The images in (b), (c) and (d) are obtained through a threshold transform of the original picture added Gaussian noise of three different noise intensities. It is clear from (b) to (d) that a suitable amount of noise reveals the hidden keyboard image but superfluous noise conceals it again.

Several methods have been presented to enhance the resonance effects for various applications, for instance, altering the shape of the potential and periodic forces of a single system [Wio & Bouzat, 1999; Gandhimathi *et al.*, 2008; Arathi *et al.*, 2011] or coupling several subsystems to obtain an enhanced collective response [Jung *et al.*, 1992; Lindner *et al.*, 1995; Morillo *et al.*, 1995; Kang & Jiang, 2009; Nicolis & Nicolis, 2017]. Specifically, the earliest investigation on the collective response can be dated back to the work based on nonlinear master equations [Jung *et al.*, 1992]. Later on, enhanced SR was applied to a nearest-neighbor coupled system of overdamped nonlinear oscillators [Lindner *et al.*, 1995]. Further on, in the mean-field coupling limit of overdamped oscillators, Gaussian white noise-induced SR of the order parameter was investigated in [Morillo *et al.*, 1995; Kang & Jiang, 2009]. Very recently, it was found in [Nicolis & Nicolis, 2017] that the presence of spatial degrees-of-freedom modifies the transition mechanism of the finite-size mean-field coupled model of overdamped oscillators. All these investigations revealed that coupling is a very effective way to enhance the effects of SR.

The prototypical system for SR is the noisy overdamped bistable oscillator model. And yet it was found that the underdamped bistable oscillator model is more effective in the early detection of mechanical defaults [Rebolledo-Herrera & Fv, 2016; López *et al.*, 2017]. This motivates the present study to investigate the resonance enhancement using coupling underdamped bistable units, specifically Duffing oscillators, thus providing some reference for relevant engineering applications. In fact, in the underdamped case, several interesting results were reported [Lindner *et al.*, 2001; Ngouongo *et al.*, 2017]. For instance, a local coupled system of nonlinear underdamped oscillators was demonstrated in [Lindner *et al.*, 2001], which exhibits multiple distinctive SR peaks by modifying the natural frequency of each subsystem.

In this paper, a diffusively coupled system of $N$ identical underdamped Duffing oscillators is considered, where the evolution of each oscillator is governed by the Langevin-type stochastic differential equation

$$\ddot{x}_i + \gamma \dot{x}_i = ax_i - bx_i^3 + \xi_i(t) + \frac{\mu}{N}\sum_{j=1}^{N}(x_j - x_i) + \varepsilon(t), \quad 1 \leq i \leq N \tag{1}$$

where the constants $a$ and $b$ are positive, so that the deterministic autonomous subsystem is bistable, $\varepsilon(t)$ represents an external periodic signal, with white Gaussian noise $\xi_i(t)$ satisfying $\langle \xi_i(t) \rangle = 0$ and $\langle \xi_i(t+\tau)\xi_j(t) \rangle = 2D\delta_{ij}(\tau)$, $i, j = 1, 2, \ldots, N$, which models the environmental fluctuations, $\gamma$ is the damping coefficient and $\mu$ is the coupling strength of the mean-field interactions among subsystems. Usually, the mean output response defined by $X(t) = \frac{1}{N}\sum_{i=1}^{N} x_i(t)$, also referred to as the order parameter, is adopted to describe the collective response of the system.

It is noted that model (1) was used to describe the dynamics of muscle contraction [Kometani & Shimizu, 1975; Desai & Zwanzig, 1978; Dawson, 1983] when the external perturbation is absent, i.e. $\varepsilon(t) = 0$. In the present paper, the focus is on the case with a sufficiently large value of *N*. In this case, with the law of large numbers, all the subsystems of model (1) can be simplified to the same evolution equation:

$$\ddot{x} + \gamma \dot{x} = (a - \mu)x - bx^3 + \mu X(t) + \xi(t) + \varepsilon(t). \tag{2}$$

Let $P(x, v, t)$ be the coordinate-velocity probability density function at time $t$ of model (2). Then, the nonlinear Fokker-Planck (FP) equation governing the probability density evolution reads [Cardiner, 1985]

$$\frac{\partial P}{\partial t} = \gamma D \frac{\partial^2 P}{\partial v^2} + \frac{\partial}{\partial v}\{[\gamma v - (a-\mu)x + bx^3 - \mu X(t) - \varepsilon(t)]P\} - v\frac{\partial P}{\partial x} \tag{3}$$

where $X(t) = \iint_{R^2} xP(x,v,t)dxdv$ with $v = \dot{x}$ and $P(x,v,t)$ obeys the natural boundary conditions. It is easy to verify that the stationary solution of Eq. (3) when $\varepsilon(t) = 0$ is

$$P_0(x,v) = Z^{-1}\exp\left\{\left(-\frac{b}{4}x^4 + \frac{1}{2}(a-\mu)x^2 + \mu X_0 x - \frac{v^2}{2}\right)\Big/D\right\} \tag{4}$$

with a normalization constant $Z$ satisfying $\iint_{R^2} P_0(x,v)dxdv = 1$. Here, the equilibrium order parameter $X_0$ is implicitly determined by $X_0 = \iint_{R^2} xP_0(x,v)dxdv$.

In Sec. 2 next, it will be proven that a Boltzmann *H*-theorem holds for the nonlinear FP

Eq. (3). Then, bifurcation diagrams of the equilibrium order parameter are depicted and analyzed. In Sec. 3, under the framework of the general linear response theory, applied to the nonlinear FP equation, an explicit relation between the order parameter and the linear susceptibility is derived. To that end, a procedure for calculating the order parameter in terms of the linear mean response is devised. Furthermore, the mono- or double-resonance behavior of the order parameter is discussed. In Sec. 4, the connection between the resonance peak and the order parameter bifurcation is revealed. Finally, conclusions are drawn in Sec. 5.

## 2. A Boltzmann-type $H$-theorem

The classical Boltzmann $H$-theorem states that the transient solution starting from any initial condition of an autonomous nonlinear evolutionary equation converges to one of the asymptotical solutions as time approaches infinity. It was developed to capture the tendency for a gas to return to an equilibrium of the evolution equation [Huang, 1963]. A Boltzmann-type $H$-theorem for the nonlinear FP equation corresponding to the mean-field coupling limit of overdamped oscillators was developed in [Shiino, 1985, 1987]. Here, the $H$-theorem is generalized from the overdamped setting to the underdamped setting when $\varepsilon(t) = 0$.

Assume that $P(x, v, t)$ is an arbitrary solution of Eq. (3) with $\varepsilon(t) = 0$ therein. Then, inspired by the work in [Shiino, 1985, 1987], define a new $H$ functional by

$$H(P(\cdot,t)) = \iint_{R^2} P(x,v,t) \ln \frac{P(x,v,t)}{Q(x,v,t)} dxdv, \tag{5}$$

where $Q(x,v,t) = \exp\left\{\left(-\frac{b}{4}x^4 + \frac{1}{2}(a-\mu)x^2 + \mu Xx - \frac{v^2}{2} - \frac{\mu X^2}{2}\right)\middle/ D\right\}$ and $X = X(t)$ is the order parameter of the unperturbed system. Comparison between $Q(x, v, t)$ and $P_0(x,v)$ yields a stationary condition,

$$\gamma D \frac{\partial^2 Q}{\partial v^2} + \frac{\partial}{\partial v}\left\{[\gamma v - (a-\mu)x + bx^3 - \mu X]Q\right\} - v \frac{\partial Q}{\partial x} = 0. \tag{6}$$

Moreover, it is easy to verify that

$$\iint_{R^2} \frac{P(x,v,t)\partial_t Q(x,v,t)}{Q(x,v,t)} dxdv = 0. \tag{7}$$

Then, based on Eqs. (6) and (7), the following Boltzmann-type $H$-theorem can be established.

**$H$-Theorem** The $H$ functional $H(P(\cdot,t))$ defined in (5) has the following properties:

(i) $H(P(\cdot,t))$ is bounded from below, namely, $H(P(\cdot,t)) > const$;

(ii) $\dfrac{dH(P(\cdot,t))}{dt} = -\iint_{R^2} \gamma D P(x,v,t) \left(\dfrac{\partial}{\partial v} \ln \dfrac{P(x,v,t)}{Q(x,v,t)}\right)^2 dxdv \leq 0$, that is, $H(P(\cdot,t))$ as a

function of time is monotonically decreasing.

Proof: (i) Let $Q_1(x,v,t)$ and $Q_2(x,v,t)$ be any time-dependent solutions of the unperturbed nonlinear FP equation. Then, $Q_2 - Q_1 \geq Q_1 \ln(Q_2/Q_1)$ since $x - 1 \geq \ln x\ (x > 0)$ and, thus,

$$\iint_{R^2} Q_1 \ln(Q_1/Q_2) dx dv \geq 0. \tag{8}$$

Here, the normalization requirements for the two probability density functions have been used. Define

$$Z(X) = \iint_{R^2} Q(x,v,t) dx dv = \iint_{R^2} \exp\left\{\left(-\frac{b}{4}x^4 + \frac{1}{2}(a-\mu)x^2 + \mu X x - \frac{v^2}{2} - \frac{\mu X^2}{2}\right)\Big/ D\right\} dx dv.$$

Substituting $Q_1(x,v,t) = P(x,v,t)$ and $Q_2(x,v,t) = Q(x,v,t)/Z(X)$ into Eq. (8) yields

$$\iint_{R^2} P \ln(P/Q) dx dv \geq -\ln Z(X). \tag{9}$$

Note that $-\dfrac{\mu}{2D}X^2 + \ln \int_R \exp\{(-\dfrac{b}{4}x^4 + \dfrac{1}{2}(a-\mu)x^2 + \mu X x)/D\} dx \leq c_0$ [Shiino, 1985] and

$\ln \int_R \exp\left\{-\dfrac{v^2}{2D}\right\} dv = \dfrac{1}{2}\ln 2\pi + \dfrac{1}{2}\ln D < c_1$. Thus, $\ln Z(X) \leq c_0 + c_1 = c$. Hence, from Eq. (9) it follows that the $H$ functional is bounded from below, namely,

$$H(P(\cdot,t)) = \iint_{R^2} P(x,v,t) \ln \frac{P(x,v,t)}{Q(x,v,t)} dx dv \geq -\ln Z(X) \geq -c. \tag{10}$$

(ii) Let $P(x,v,t)$ be a solution of the nonlinear FP Eq. (3) with $\varepsilon(t) = 0$ therein. Then, taking the derivative of $H(P(\cdot,t))$ with respect to time $t$, one obtains

$$\frac{dH(P(\cdot,t))}{dt} = \iint_{R^2} \dot{P} \ln \frac{P}{Q} dx dv + \iint_{R^2} P\left(\frac{\dot{P}}{P} - \frac{\dot{Q}}{Q}\right) dx dv$$

$$= \iint_{R^2} \left\{\gamma D \frac{\partial^2 P}{\partial v^2} + \frac{\partial}{\partial v}\left[(\gamma v - (a-\mu)x + bx^3 - \mu X)P\right] - v\frac{\partial P}{\partial x}\right\} \ln \frac{P}{Q} dx dv - \iint_{R^2} \frac{P\dot{Q}}{Q} dx dv$$

$$= \iint_{R^2} P\left\{\gamma D \frac{\partial}{\partial v} - [\gamma v - (a-\mu)x + bx^3 - \mu X]\right\} \frac{\partial/\partial v(P/Q)}{P/Q} dx dv$$

$$+ \iint_{R^2} Pv \frac{\partial/\partial x(P/Q)}{P/Q} dx dv - \iint_{R^2} \frac{P\dot{Q}}{Q} dx dv$$

$$= \iint_{R^2} \frac{P}{P/Q}\left\{\gamma D \frac{\partial}{\partial v} - [\gamma v - (a-\mu)x + bx^3 - \mu X]\right\} \frac{\partial}{\partial v}(P/Q) dx dv - \iint_{R^2} P\gamma D \left(\frac{\partial/\partial v(P/Q)}{P/Q}\right)^2 dx dv$$

$$+ \iint_{R^2} Pv \frac{\partial/\partial x (P/Q)}{P/Q} dxdv - \iint_{R^2} \frac{P\dot{Q}}{Q} dxdv$$

$$= \iint_{R^2} \frac{P}{Q} \left\{ \gamma D \frac{\partial^2}{\partial v^2} + \frac{\partial}{\partial v} \{ [\gamma v - (a-\mu)x + bx^3 - \mu X] \} - v \frac{\partial}{\partial x} \right\} Q dxdv$$

$$- \iint_{R^2} P\gamma D \left( \frac{\partial}{\partial v} \ln \frac{P}{Q} \right)^2 dxdv - \iint_{R^2} \frac{P\dot{Q}}{Q} dxdv. \tag{11}$$

According to Eqs. (6) and (7), the first and third terms on the right-hand side of Eq. (11) vanish. Therefore, one finally obtains

$$\frac{dH(P(\cdot,t))}{dt} = -\iint_{R^2} \gamma PD \left( \frac{\partial}{\partial v} \ln \frac{P}{Q} \right)^2 dxdv \leq 0. \tag{12}$$

The above $H$-theorem states that the $H$ functional decreases monotonically in time due to property (ii) but cannot decrease indefinitely due to property (i). As a result, $dH(P(\cdot,t))/dt$ must vanish as time approaches infinity. Thus, from Eq. (12), one obtains $\lim_{t\to\infty} \frac{\partial}{\partial v} \ln \frac{P}{Q} = 0$, which implies that the ratio $P/Q$ is independent of $v$ as $t \to \infty$.

Let $f(x) = \lim_{t\to\infty} \frac{P}{Q}$. Then, $\lim_{t\to\infty} P(x,v,t) = f(x) \lim_{t\to\infty} Q(x,v,t)$ is an asymptotic solution of Eq. (3) in the case of $\varepsilon(t) = 0$. Substituting this asymptotic solution into Eq. (3) and noticing the definition of $Q(x, v, t)$, one obtains $\partial f(x)/\partial x = 0$ and $\lim_{t\to\infty} \partial X/\partial t = 0$. It thus follows that $f(x) = const$ and $\lim_{t\to\infty} X = const$. Hence, $Q(x,v,\infty)$ should be coincident with $P_0(x,v)$ up to a normalization constant, implying that any time-dependent solution of the unperturbed FP equation (3) converges to one of its stationary probability density functions as time tends to infinity.

Since any transient solution converges to some stationary solution in the large time limit, one can investigate the bifurcation behavior of the nonlinear Eq. (3) in the absence of $\varepsilon(t)$ by means of the implicit equation $X_0 = \iint_{R^2} x P_0(x,v) dxdv$. Fig.2 shows the bifurcation diagram of the order parameter versus the coupling strength and the noise intensity, respectively. As can be seen from the figure, for a given $\mu$, there exists a critical noise intensity $D_c$ such that the equilibrium order parameter takes two stable states when $D < D_c$ but only one stable state when $D > D_c$. For a fixed $D$, there also exists a critical coupling strength $\mu_c$ such that the equilibrium order parameter undergoes a (pitchfork) bifurcation as the value of $\mu$ increases.

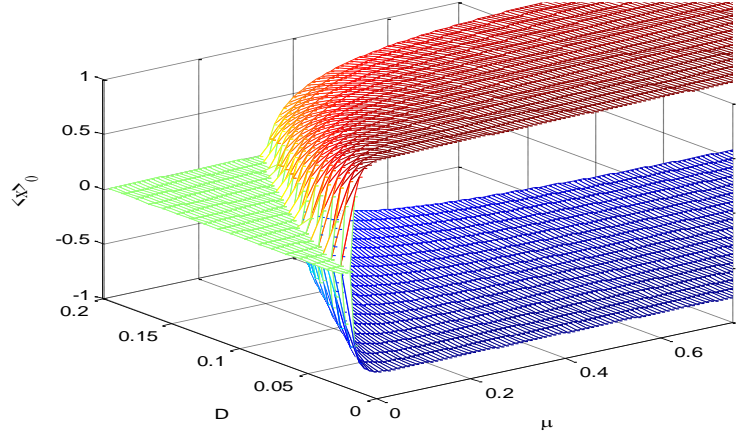

Fig.2. Bifurcation diagram of the stable equilibrium order parameter versus the coupling strength and the noise intensity. The involving system parameters are $\gamma = 0.4$, $a = 1.0$ and $b = 1.0$, respectively.

## 3. Linear response analysis

The classical Floquet theory applicable to time-varying linear systems cannot guarantee that the asymptotic solution of the nonlinear FP equation is periodic in time [Morillo *et al.*, 1995]. Nevertheless, due to the above-derived *H*-theorem, any time-dependent solution of the unperturbed version of nonlinear FP Eq. (3) can converge to one of the stationary solutions. Therefore, one can employ the technique of perturbation analysis [Risken, 1989] to explore the long-term dynamics of system (3), provided that the external input $\varepsilon(t)$ is sufficiently weak.

### *3.1 Linear susceptibility versus order parameter*

Within the linear response range of $|\varepsilon(t)| \ll 1$, the long-time probability density function can be formulated into the form of the perturbation expansion

$$P_{as}(x,v,t) = P_0(x,v) + P_1(x,v,t), \tag{13}$$

where $P_0(x,v)$ is the equilibrium solution of Eq. (3) when $\varepsilon(t) = 0$ and $P_1(x,v,t)$ is to be determined. The long-time order parameter can also be split into $X_{as}(t) = X_0 + X_1(t)$. Substituting Eq. (13) into Eq. (3) gives

$$\frac{\partial P_1(x,v,t)}{\partial t} = L_{FP}(x,v)P_1(x,v,t) + L_{ext}(x,v,t)P_0(x,v), \tag{14}$$

where $L_{FP}(x,v) = -\frac{\partial}{\partial x}v + \frac{\partial}{\partial v}(\gamma v - (a-\mu)x + bx^3 - \mu X_0) + D\gamma \frac{\partial^2}{\partial v^2}$ is the nonlinear

unperturbed FP operator and $L_{ext}(x,v,t) = -[\mu X_1(t) + \varepsilon(t)]\frac{\partial}{\partial v}$ is the perturbation operator from the external field. Then, the asymptotic solution of Eq. (14) can be formally expressed as

$$P_1(x,v,t) = \int_{-\infty}^{t} e^{L_{FP}(x,v)(t-t')} L_{ext}(x,v,t') P_0(x,v) dt'. \tag{15}$$

The deviation of the expectation of $X(t)$ with respect to the equilibrium value yields

$$\begin{aligned} X_1(t) &= \iint_{R^2} x P_1(x,v,t) dx dv \\ &= \int_{-\infty}^{t} \iint_{R^2} x e^{L_{FP}(x,v)(t-t')} [\mu X_1(t') + \varepsilon(t')] \frac{v}{D} P_0(x,v) dx dv dt' \\ &= \int_{-\infty}^{t} R_{x,v}(t-t')[\mu X_1(t') + \varepsilon(t')] dt', \end{aligned} \tag{16}$$

where

$$R_{x,v}(t) = \begin{cases} \frac{1}{D} \iint_{R^2} xv e^{L_{FP}(x,v)t} P_0(x,v) dx dv, & t \geq 0; \\ 0, & t < 0. \end{cases}$$

In the above derivation of Eq. (16), the identity $\frac{\partial P_0(x,v)}{\partial v} = -\frac{v}{D} P_0(x,v)$ has been used.

When $t \geq t_0 = 0$, one has

$$P(x,v,t|x_0,v_0,t_0) = e^{L_{FP}(x,v)(t-t_0)} \delta(x-x_0)\delta(v-v_0)$$

and

$$\frac{d}{d\bar{t}} P(x,v,t|x_0,v_0,\bar{t})\big|_{\bar{t}=t_0} = -e^{L_{FP}(x,v)(t-t_0)} L_{FP}(x,v)\delta(x-x_0)\delta(v-v_0).$$

Since $L_{FP}(x,v)P_0(x,v) = 0$, a direct calculation yields

$$\begin{aligned} R_{x,v}(t) &= \frac{1}{D} \int_{R^4} x P(x,v,t|x_0,v_0,0) P_0(x_0,v_0) v_0 dx dv dx_0 dv_0 \\ &= \frac{1}{D} \int_{R^4} \frac{d}{d\bar{t}} \{x P(x,v,t|x(\bar{t}),v(\bar{t}),\bar{t}) P_0(\bar{x}(t),\bar{v}(t)) x(\bar{t})\}_{\bar{t}=t_0} dx dv dx_0 dv_0 \\ &\quad - \frac{1}{D} \int_{R^4} x \frac{d}{d\bar{t}} [P(x,v,t|x(\bar{t}),v(\bar{t}),\bar{t}) P_0(x(\bar{t}),v(\bar{t}))]\big|_{\bar{t}=t_0} x_0 dx dv dx_0 dv_0 \\ &= -\frac{1}{D} \frac{d}{dt} K_{xx}(t) + \frac{1}{D} \iint_{R^2} x e^{L_{FP}(x,v)t} [L_{FP}(x,v) P_0(x,v)] x dx dv \\ &= -\frac{1}{D} \frac{d}{dt} K_{xx}(t) \end{aligned}$$

where $K_{xx}(t) = \langle x(t)x(0) \rangle$ is the stationary autocorrelation function for the unperturbed system. Then, the above equation leads to the well-known Klein-Kramers relation [Risken, 1989]

$$R_{x,v}(t) = \begin{cases} -\dfrac{1}{D}\dfrac{d}{dt}K_{xx}(t), & t \geq 0; \\ 0, & t < 0. \end{cases}$$

Let $\tilde{R}(\omega) = \int_{-\infty}^{+\infty} R_{x,v}(t)e^{-i\omega t}dt = \dfrac{1}{D}[K_{xx}(0) - i\omega \int_{0}^{+\infty} K_{xx}(t)e^{-i\omega t}dt]$. Then, from Eq. (16) one has $\tilde{X}_1(\omega) = \left(\mu \tilde{X}_1(\omega) + \tilde{\varepsilon}(\omega)\right)\tilde{R}(\omega)$, and thus the deviation of the order parameter caused by the external perturbation, in the frequency domain, reads

$$\tilde{X}_1(\omega) = \chi(\omega)\tilde{\varepsilon}(\omega). \tag{17}$$

Here, $\chi(\omega) = \dfrac{\tilde{R}(\omega)}{1 - \mu \tilde{R}(\omega)}$ satisfying $\chi(-\omega) = \chi^*(\omega)$ is the so-called linear dynamic susceptibility of the order parameter [Morillo *et al.*, 1995]. To this end, by letting $\varepsilon(t) = \varepsilon_0 \cos \Omega t$, it immediately follows from Eq. (17) that

$$\tilde{X}_1(\omega) = \varepsilon_0 \pi \chi(\omega)[\delta(\omega+\Omega) + \delta(\omega-\Omega)]$$

and $X_1(t) = \varepsilon_0 \text{Re}[\chi(\Omega)e^{i\Omega t}]$.

*3.2 Calculation of linear susceptibility*

Next, take $\varepsilon(t) = \varepsilon_0 \cos \Omega t$ with $\varepsilon_0 \ll 1$ to calculate the linear susceptibility. The main idea underlying the calculating procedure is to combine the method of weighed series expansion [Evstigneev *et al.*, 2002; Kang *et al.*, 2003; Kang, 2011; Liu & Kang, 2018] and the harmonic balance method [Lim *et al.*, 2001], where the weighting function is taken as the stationary solution of the unperturbed version of Eq. (3) so as to satisfy the natural boundary conditions. This method was applied to linear FP equations before [Evstigneev *et al.*, 2002; Kang *et al.*, 2003; Kang, 2011; Liu & Kang, 2018], but it is applied to the nonlinear FP equation here for the first time.

Under the framework of the linear response theory, one can seek the asymptotic solution of Eq. (13) such that

$$P_1(x,v,t) = \varepsilon_0(P_{11}(x,v)\cos\Omega t + P_{12}(x,v)\sin\Omega t), \tag{18}$$

where $P_{11}(x,v)$ and $P_{12}(x,v)$ are unknown functions satisfying

$$\iint_{R^2} P_{11}(x,v)dxdv = \iint_{R^2} P_{12}(x,v)dxdv = 0.$$

due to the normalization requirement as a probability.

Substituting Eq. (18) into Eq. (3) and using the orthogonality of trigonometric functions, one obtains

$$-\Omega P_{11}(x,v) = \gamma D \frac{\partial^2}{\partial v^2} P_{12}(x,v) + \frac{\partial}{\partial v}\{[\gamma v - (a-\mu)x + bx^3 - \mu\langle x\rangle_0]P_{12}(x,v)\}$$

$$+ \frac{\partial}{\partial v}[-\mu\langle x\rangle_{12}P_0(x,v)] - v\frac{\partial}{\partial x}P_{12}(x,v) \qquad (19)$$

and

$$\Omega P_{12}(x,v) = \gamma D \frac{\partial^2}{\partial v^2} P_{11}(x,v) + \frac{\partial}{\partial v}\{[\gamma v - (a-\mu)x + bx^3 - \mu\langle x\rangle_0]P_{11}(x,v)\}$$

$$+ \frac{\partial}{\partial v}\{[-\mu\langle x\rangle_{11} - 1]P_0(x,v)\} - v\frac{\partial}{\partial x}P_{11}(x,v). \qquad (20)$$

Here, $X_{as}(t) = \langle x\rangle_0 + \varepsilon_0(\langle x\rangle_{11}\cos\Omega t + \langle x\rangle_{12}\sin\Omega t)$ with $\langle\cdots\rangle_i$ and $\langle\cdots\rangle_{1i}$ stands for the integral with respect to the probability density $P_i(x,v)$ and $P_{1i}(x,v)$, respectively.

Let $P_1(x,v) = P_{11}(x,v) + iP_{12}(x,v)$. Then, Eqs. (19) and (20) can be combined into one equation as

$$-i\Omega P_1(x,v) = \gamma D \frac{\partial^2}{\partial v^2} P_1(x,v) + \frac{\partial}{\partial v}\{[\gamma v - (a-\mu)x + bx^3 - \mu\langle x\rangle_0]P_1(x,v)\}$$

$$+ \frac{\partial}{\partial v}\{[-\mu\langle x\rangle_1 - 1]P_0(x,v)\} - v\frac{\partial}{\partial x}P_1(x,v). \qquad (21)$$

Let $F(x,v)$ be an arbitrary function of coordinates $x$ and $v$, and assume that the corresponding time-dependent moment $\langle F(x,v)\rangle(t) = \iint_{R^2} F(x,v)P(x,v,t)dxdv$ exists. Then, multiplying the right-hand and the left-hand sides of Eq. (21) by $F(x,v)$ and then integrating the result over the $x-v$ plane leads to

$$-i\Omega\langle F\rangle_1 = \gamma D\left\langle\frac{\partial^2 F}{\partial v^2}\right\rangle_1 - \left\langle\{[\gamma v - (a-\mu)x + bx^3 - \mu\langle x\rangle_0]\}\frac{\partial F}{\partial v}\right\rangle_1 + \left\langle v\frac{\partial F}{\partial x}\right\rangle_1 + \left\langle\frac{\partial F}{\partial v}\right\rangle_0 + \left\langle\mu\langle x\rangle_1 \frac{\partial F}{\partial v}\right\rangle_0.$$

Rewriting $P_1(x,v) = P_0(x,v)p_1(x,v)$ and noting that

$$\frac{\partial P_0(x,v)}{\partial x} = \frac{1}{D}[-bx^3 + (a-\mu)x + \mu\langle x\rangle_0]P_0(x,v), \quad \frac{\partial P_0(x,v)}{\partial v} = -\frac{v}{D}P_0(x,v),$$

one obtains

$$-i\Omega\langle Fp_1\rangle_0 = -D\gamma\left\langle\frac{\partial F}{\partial v}\frac{\partial p_1}{\partial v}\right\rangle_0 + D\left\langle\frac{\partial F}{\partial x}\frac{\partial p_1}{\partial v} - \frac{\partial F}{\partial v}\frac{\partial p_1}{\partial x}\right\rangle_0 + \left\langle\frac{\partial F}{\partial v}\right\rangle_0 + \mu\langle xp_1\rangle_0\left\langle\frac{\partial F}{\partial v}\right\rangle_0. \qquad (22)$$

In order to find $p_1(x,v)$, decompose $p_1(x,v)$ into

$$p_1(x,v) = \sum_{j=0}^{\infty}\sum_{k=0}^{\infty} c_{j,k} H_j(v/v_t)x^k, \qquad (23)$$

where $H_j$ stands for the $j$th Hermite polynomial and $v_t = \sqrt{2D}$.

Here, it is remarked that the decomposition with respect to the velocity $v$ is the same as that in [Evstigneev *et al.*, 2002; Liu & Kang, 2018], but a bit different from that in [Kang *et al.*, 2003; Kang, 2011]. Nevertheless, numerical experiment shows that the two forms of decomposition have the same convergence speed.

Now, substituting expansion (23) into Eq. (22) with $F(x,v) = H_s(v/v_t)x^l$, $s = 0, 1, \ldots$, one obtains an infinite set of linear algebraic equations:

$$\sum_{s=0}^{\infty}\sum_{k=0}^{\infty} c_{s,k}(-i\Omega + \gamma s)\langle x^{k+l}\rangle_0 = \sum_{s=0}^{\infty}\sum_{k=0}^{\infty} c_{s+1,k}\sqrt{2D}(s+1)l\langle x^{k+l-1}\rangle_0 - \sum_{s=1}^{\infty}\sum_{k=0}^{\infty} c_{s-1,k}\sqrt{D/2}k\langle x^{k+l-1}\rangle_0$$

$$+ \mu\left[\sum_{k=0}^{\infty} c_{0,k}\langle x^{k+1}\rangle_0\right]\cdot\left[\langle x^l\rangle_0/\sqrt{2D}\right] + \langle x^l\rangle_0/\sqrt{2D}, \qquad (24)$$

where some properties of Hermite polynomials have been applied. For numerical calculation, truncate $j$ and $k$ in Eq. (23) and then substitute it into Eq. (24). Then, a block-tridiagonal system for the unknown coefficient vectors $c_s = (c_{s,0}, c_{s,1}, \ldots, c_{s,K})^T$ is obtained as follows:

$$\begin{bmatrix} A_0 & B_0 & & & \\ C_1 & A_1 & B_1 & & \\ & C_2 & A_2 & \ddots & \\ & & \ddots & \ddots & B_{J-1} \\ & & & C_J & A_J \end{bmatrix}\begin{bmatrix} c_0 \\ c_1 \\ c_2 \\ \vdots \\ c_J \end{bmatrix} = \begin{bmatrix} 0 \\ d_1 \\ 0 \\ \vdots \\ 0 \end{bmatrix} \qquad (25)$$

with entries of the involving matrices given by

$$(A_s)_{l,k} = (\gamma s - i\Omega)\langle x^{l+k}\rangle_0, \quad (B_s)_{l,k} = -\sqrt{2D}(s+1)l\langle x^{l+k-1}\rangle_0,$$

$$(C_1)_{l,k} = \sqrt{D/2}k\langle x^{l+k-1}\rangle_0 - \mu\frac{1}{\sqrt{2D}}\langle x^{k+1}\rangle_0\langle x^l\rangle_0,$$

$$(C_s)_{l,k} = \sqrt{D/2}k\langle x^{l+k-1}\rangle_0 \quad (s \geq 2),$$

and the components of the right-hand vector given by $(d_1)_l = \langle x^l\rangle_0/\sqrt{2D}$. When $s = 0$ and $l = 0$, the normalization of the probability is satisfied naturally.

The Gaussian block-elimination method can then be used to solve system (25). With the solution so obtained, the linear dynamical susceptibility $\chi(\Omega) = \iint_{R^2} xP_1(x,v)dxdv$ can be obtained, as

$$\chi(\Omega) = \langle xp_1\rangle_0 = \sum_{k=0}^{K} c_{0,k}\langle x^{k+1}\rangle_0. \qquad (26)$$

The numerical experiment shows that the method of moments converges quite fast as the truncation orders increase. In the following calculations, fix $K=10, J=10$ because larger values of $K$ and $J$ do not change the results by more than 0.5% in the parameter range considered. To confirm the accuracy of the method of moments, the Euler-Maruyama method [Fox *et al.*, 1988; Higham, 2001] is also adopted for simulations. As shown in Fig. 3, there is a good agreement between the theoretical curve and the simulated curve, validating the method of moments with the given parameters. In fact, more comparisons have been carried out, showing a good agreement between the two methods for $D \geq 0.03$. Corresponding theoretical results include the one shown in Fig. 4, where only the theoretic results with noise intensity are displayed under this constraint.

With the linear dynamical susceptibility available, one can use the spectral amplification factor defined by $|\chi(\Omega)|^2$ [Jung, 1993; Gammaitoni *et al.*, 1998] to measure the resonance behavior in the order parameter of system (2). It was found [Alfonsi *et al.*, 2000; Kang *et al.*, 2003] that, as a result of the coexistence of three types of motions (i.e., intrawell oscillation, interwell random transition and over-barrier vibration), the dependence of the spectral amplification factor on the noise intensity in a single underdamped bistable Duffing oscillator might exhibit mono- or double-peak resonance structures.

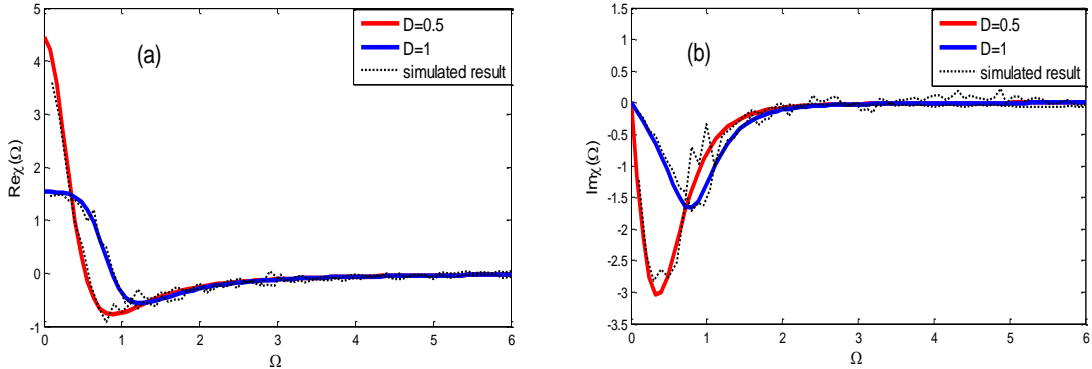

Fig.3. Dependence of the linear dynamical susceptibility (real part (a) and imaginary part (b)) on driving frequency. The system parameters are taken as $a=1.0$, $b=1.0$, $\gamma=0.4$ and $\mu=0.6$.

Here, the question is how the mean-filed coupling will affect these structures. As can be seen from Fig. 4, the order parameter can once again exhibit double-peak resonance behavior at some driving frequencies. Interestingly, the effects of the mono-peak resonance and the double-peak resonance are both enhanced. This observation enables one to generalize the existing results of coupling enhancing SR [Jung *et al.*, 1992; Lindner *et al.*, 1995; Morillo *et al.*, 1995; Kang & Jiang, 2009; Nicolis & Nicolis, 2017], from the overdamped setting to the underdamped setting. More interestingly, by tuning the mean-coupling strength with a fixed noise intensity, a weird parameter-induced double resonance can be observed. This implies that tuning parameter paves a new way to induce double-peak SR.

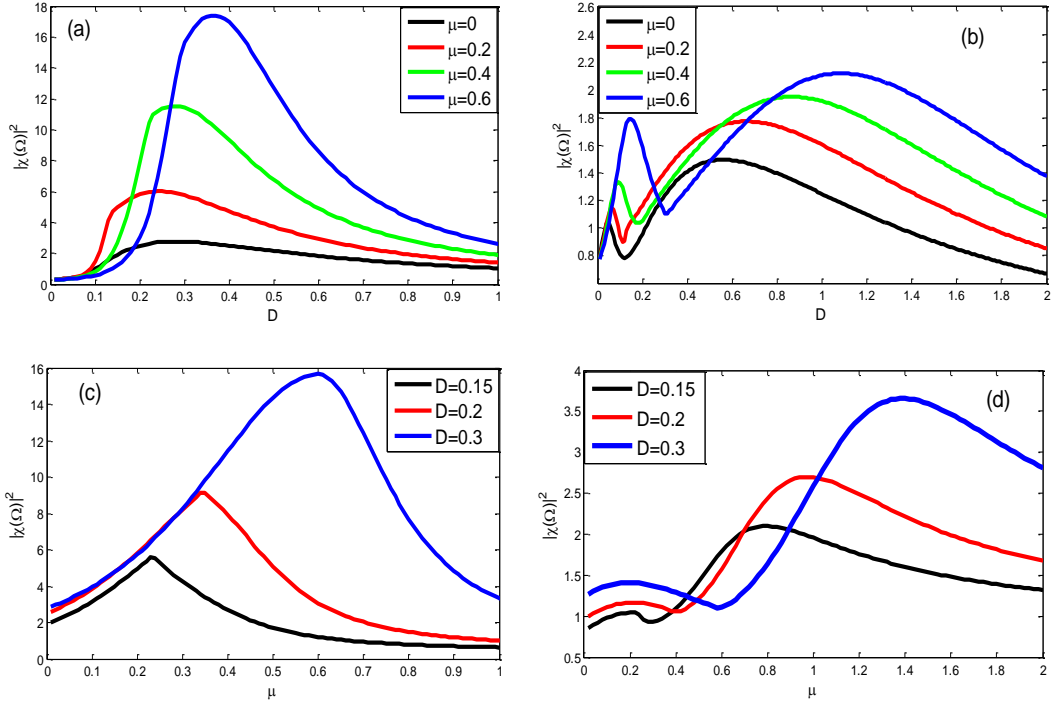

Fig.4. Dependence of the spectral amplification factor on the noise intensity (a), (b) and on the coupling strength (c), (d). The other parameters are taken as $a=1.0, b=1.0, \gamma=0.4$ and $\Omega=0.1\pi$ in (a), (c), $\Omega=0.3\pi$ in (b), (d).

### 3.3 Response to general weak periodic signals

Consider the question of how to obtain the long-time response of the order parameter to general weak periodic signals using the principle of superposition. Two examples are discussed.

First, consider the case where the external perturbation is a weak periodic square wave signal of the form $\varepsilon(t)=(-1)^{k-1}\varepsilon_0, t\in\left[\dfrac{\pi(k-1)}{\Omega},\dfrac{\pi k}{\Omega}\right], k=1,2,\ldots$ with $\varepsilon_0\ll 1$. Noting that the Fourier series for this signal reads $\varepsilon(t)=\sum_{k=1}^{\infty}M\sin((2k-1)\Omega t)$ with $M=\dfrac{4\varepsilon_0}{\pi(2k-1)}$, one immediately obtains the following Fourier transform:

$$\tilde{\varepsilon}(\omega)=\sum_{k=1}^{\infty}Mi\pi[\delta(\omega+(2k-1)\Omega)-\delta(\omega-(2k-1)\Omega)].$$

Substituting this equation into Eq. (17) yields

$$\tilde{X}_1(\omega)=\chi(\omega)\sum_{k=1}^{\infty}Mi\pi[\delta(\omega+(2k-1)\Omega)-\delta(\omega-(2k-1)\Omega)].$$

Thus, an operation of inverse Fourier transform gives

$$X_1(t)=\dfrac{i}{2}\sum_{k=1}^{\infty}M[\chi((2k-1)\Omega)\exp(i(2k-1)\Omega t)-\chi(-(2k-1)\Omega)\exp(-i(2k-1)\Omega t)].$$

Hence, the long-time order parameter within linear response range (as shown in Fig. 5) reads

$$X_{as}(t) = X_0 + \sum_{k=1}^{\infty} M \, \text{Im}[\chi((2k-1)\Omega)\exp(i(2k-1)\Omega t)].$$

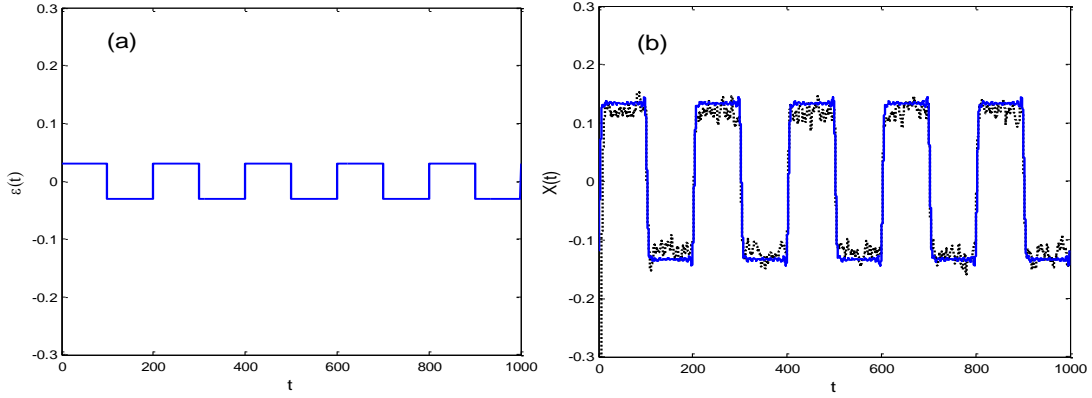

Fig.5. The periodic telegraph signal (a) and long-time order parameter (b): direct simulation (black dotted) and theoretical method (blue solid) in the first example with parameters $a=1.0$, $b=1.0$, $\mu=0.6, D=0.5, \gamma=0.4, \varepsilon_0 = 0.03$ and $\Omega=0.01\pi$.

Second, consider the case where the external perturbation is a periodic envelope signal generated from unilateral attenuation impulse in gear fault [Zhang *et al.*, 2017], in the form of $\varepsilon(t) = \varepsilon_0 \exp(-d \bmod(t,T_0))$, where $d$ is the attenuation rate and $T_0$ is the driving period. In this case, the corresponding Fourier transform is

$$\tilde{\varepsilon}(\omega) = a_0\pi\delta(\omega) + \pi\sum_{k=1}^{\infty} a_k[\delta(\omega+k\Omega)+\delta(\omega-k\Omega)] + i\pi\sum_{k=1}^{\infty} b_k[\delta(\omega+k\Omega)-\delta(\omega-k\Omega)], \quad (27)$$

where $a_k = \dfrac{2\varepsilon_0(1-\exp(-dT_0))d}{T_0(d^2+(k\Omega)^2)}$, $b_k = \dfrac{2\varepsilon_0(1-\exp(-dT_0))(k\Omega)}{T_0(d^2+(k\Omega)^2)}$ for $k=0,1,2,\ldots$.

Substitution of Eq. (27) into Eq. (17) yields

$$\tilde{X}_1(\omega) = \chi(\omega)a_0\pi\delta(\omega) + \pi\chi(\omega)\sum_{k=1}^{\infty} a_k[\delta(\omega+k\Omega)+\delta(\omega-k\Omega)] + i\pi\chi(\omega)\sum_{k=1}^{\infty} b_k[\delta(\omega+k\Omega)-\delta(\omega-k\Omega)].$$

Thus,

$$X_1(t) = \frac{a_0}{2}\chi(0) + \sum_{k=1}^{\infty} a_k \, \text{Re}[\chi(k\Omega)\exp(ik\Omega t)] + \sum_{k=1}^{\infty} b_k \, \text{Im}[\chi(k\Omega)\exp(ik\Omega t)].$$

Within the linear response range, the long-time order parameter reads $X_{as}(t) = X_0 + X_1(t)$, as shown in Fig. 6.

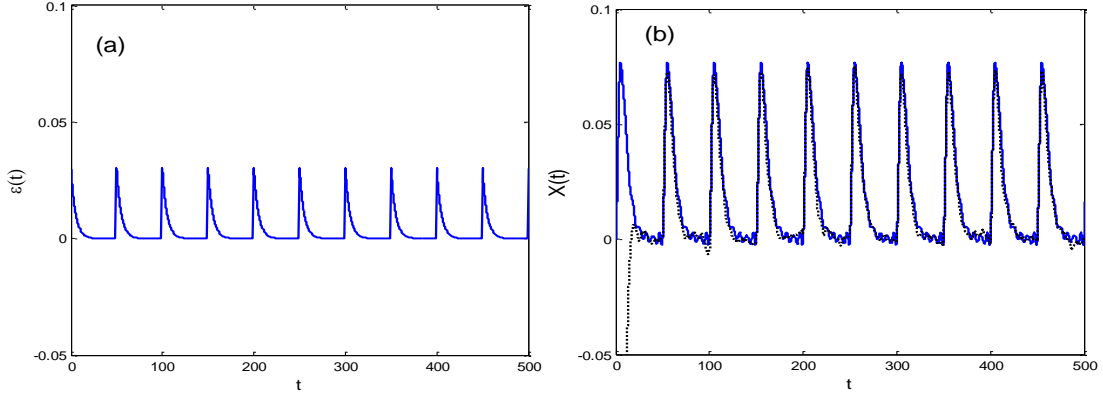

Fig.6. The envelope signal (a) and the long-time order parameter (b): direct simulation (black dotted) and theoretical method (blue solid) in the second example with parameters $a = 1.0$, $b = 1.0, \mu = 0.6$,

$$D = 0.45, \gamma = 0.4, \ \varepsilon_0 = 0.03, \ \Omega = 2\pi/T_0 = 0.04\pi \text{ and } d = 0.2.$$

Recall that Fig. 3 clearly shows that both the real part and the imaginary part of the linear dynamical susceptibility are generally more sensitive to lower frequency signals, and thus the low frequency signals should be easier to be amplified by the proposed mean-field model at a suitable noise level. In fact, as shown in Fig. 5 and Fig. 6, both low-frequency periodic signals are amplified by more than two times at the prescribed noise level.

## 4. Stochastic resonance and bifurcation

As shown in Fig.2, when the control parameter such as the noise intensity or the coupling strength passes through a critical value, the equilibrium order parameter of the mean-field coupled system (2) with the absence of external periodic signals will change from a trivial phase into a nontrivial phase. Correspondingly, the effective potential $U(x) = \frac{(\mu - a)x^2}{2} + \frac{bx^4}{4} - \mu X_0 x$ will lose its symmetry in shape at the critical point, and this will certainly affect the resonance. Thus, it is natural to question whether there is a relationship between the SR behavior of the collective response and the bifurcation of the order parameter. For simplicity, only the connection of the mono-SR peak and the order parameter bifurcation is considered here. In the following numerical simulation, time step is 0.01.

The dependence of the unperturbed order parameter and the spectral amplification factor on the noise intensity for different coupling strengths is shown in Figs.7 (a) and (b), respectively. It is easy to see that, as the mean-field coupling strength is enhanced, although the critical noise intensity at the bifurcation point and the optimal noise intensity, where the spectral amplification factor attains peak value, shifts towards a higher noise level, the SR peak grows higher and higher quickly, just like what happened in the overdamped case [Morillo *et al.*, 1995; Kang & Jiang, 2009]. This observation implies that one can firmly improve the effect of the spectral-amplification-factor measured resonance by strengthening the coupling interactions among these underdamped oscillators. Moreover, by a more careful examination, it can be seen that the critical noise intensity (denoted by the vertical dash line in Fig. 7(a)) is almost the same as the optimal noise intensity (denoted by the vertical dash line in Fig. 7(b)) for a given coupling strength. It should be emphasized that this observation is useful since it is instructive for mechanical engineers to approximate the optimal noise intensity by the critical noise intensity in their fault detection applications [Vania & Pennacchi, 2004; Leng *et al.*, 2006; Qiao *et al.*, 2017; Zhang *et al.*, 2017; Ma *et al.*, 2018].

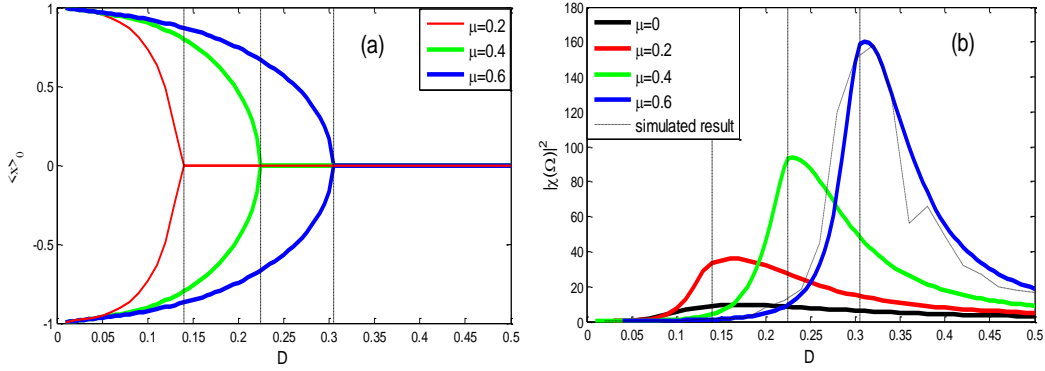

Fig.7. Bifurcation diagrams of the equilibrium order parameter versus the noise intensity (a) and dependence of the spectral amplification factor on the noise intensity (b). The system parameters are taken as $a = 1.0, b = 1.0, \Omega = 0.1$ and $\gamma = 0.4$. The vertical dash lines mark the locations of the bifurcation points and the resonant peaks. By checking the horizontal coordinates of the intersection points of the vertical dash lines with the level axis, a corresponding relation between the pitchfork bifurcation and SR is clear.

The dependence of the unperturbed order parameter and the spectral amplification factor on the coupling strength for different noise intensities is shown in Figs. 8 (a) and (b), respectively. It can be observed that the critical coupling strength at the bifurcation point becomes larger as the noise intensity becomes higher. This observation sufficiently demonstrates that the bifurcation is closely related to noise. Meanwhile, one can see that, as the noise intensity becomes higher, not only the optimal coupling strength at the resonance peak becomes larger, but also the peak height of the SR curve becomes higher. This undoubtedly demonstrates the counterintuitive role of noise in the way just like SR in a single underdamped system. Further, by checking the vertical dash lines, one can find again that the optimal coupling strength of SR is in coincidence with the critical coupling strength at the bifurcation point. Hence, this might offer a good suggestion to the mechanical engineers to approximate the optimal coupling strength by the critical coupling strength if they adopt a coupling strength parameter induced SR in their applications.

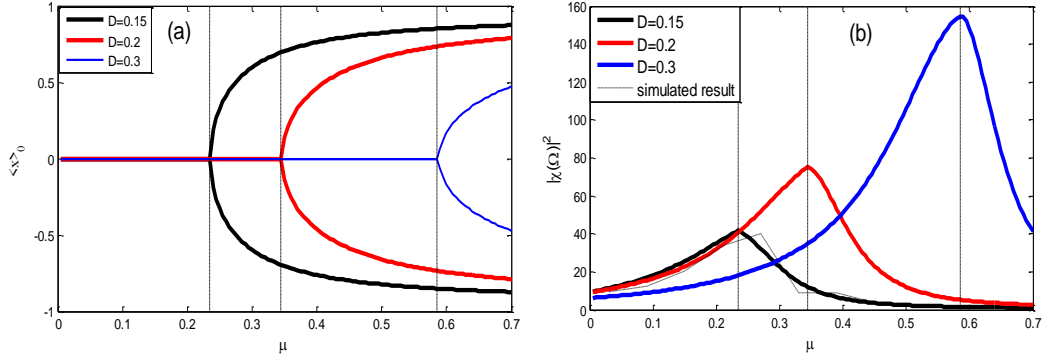

Fig.8. Bifurcation diagram of the equilibrium order parameter versus the coupling strength (a) and dependence of the spectral amplification factor on the coupling strength (b). The system parameters are taken as $a = 1.0, b = 1.0$, $\Omega = 0.1$, and $\gamma = 0.4$. Again, by checking the intersection points of the vertical dash lines with the horizontal axis, from where a corresponding relation between the horizontal coordinates is clear.

With the damping coefficient as a tunable parameter, the underdamped bistable Duffing oscillator has a more practical advantage comparing to the overdamped oscillator, since the damping coefficient might induce SR [Evstigneev *et al.*, 2002; Liu & Kang, 2018]. In order to check the effect of the mean-field coupling on the dissipation-induced SR, the dependence of the spectral amplification factor on the damping coefficient is shown in Fig. 9. From Fig.9 (a), one can see that, for fixed $\mu$, the resonance effect is not monotonically changing with the noise intensity, instead there exists an optimal noise intensity $D^{opt} = 0.16$ such that the highest resonance peak can be achieved. Analogously, for fixed $D$, from Fig. 9(b) one can see that there exists an optimal coupling strength $\mu^{opt} = 0.23$, where the highest resonance peak arrives. Moreover, a more careful scrutiny of Figs. 7-8 with Fig. 9 discloses that, for the dissipation-induced SR, the optimal noise intensity $D^{opt}$ is close to the critical noise intensity $D_c = 0.14$ of the noise intensity induced bifurcation, and the optimal coupling strength $\mu^{opt}$ takes the same value with the critical coupling strength $\mu_c$ of the coupling strength induced bifurcation. This observation further confirms the connection of the bifurcation and the SR.

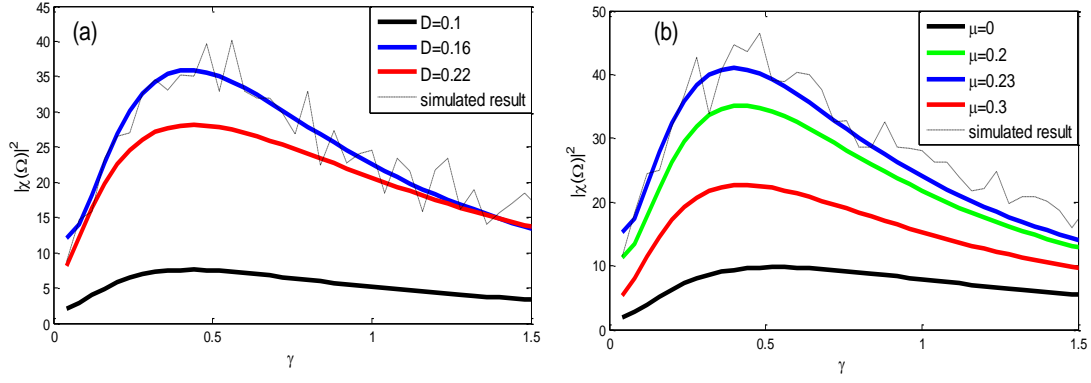

Fig.9. Dependence of the spectral amplification factor on the damping coefficient. The other parameters are taken as $a=1.0$, $b=1.0$, $\Omega=0.1$, $\mu=0.2$ in (a) and $D=0.15$ in (b). From these cures, one can see that the peak of the dissipation-induced SR is neither a monatomic function of the noise intensity nor the coupling strength. In fact, the peak is related to the critical values of the bifurcation of the order parameter showed in Figs. 7-8.

## 5. Conclusions

The long-term mean-field dynamics of coupled underdamped Duffing oscillators driven by an external periodic signal with Gaussian noise has been investigated. A Boltzmann-type $H$-theorem has been established for the associated nonlinear Fokker-Planck equation to ensure the long-time stationary state of the system. Based on the general framework of linear response theory and in terms of linear dynamical susceptibility, a key order parameter has been introduced and then the effect of the system parameters such as noise intensity and coupling strength on SR and bifurcation dynamics of the collective response has been investigated, revealing the connection between the optimal stochastic resonance parameter and the critical bifurcation parameter. The new results can be regarded as an extension and complement from the coupled model of overdamped bistable oscillators to the underdamped setting.

The main research findings can be summarized as follows. Firstly, for both cases of mono-peak and double-peak resonance curves, the mean-field coupling is always beneficial for the resonance enhancement. Secondly, for the mono-peak stochastic resonance, the optimal noise intensity is closely related to the critical noise intensity of the order parameter bifurcation; similarly, for the coupling-strength-induced mono-peak stochastic resonance, the optimal coupling strength is the critical coupling strength. Thirdly, for the dissipation-induced stochastic resonance, there exists an optimal noise intensity or an optimal coupling strength that maximizes the resonance effect. These new observations offer a good reference for practical applications, including weak signal detection and early mechanical fault diagnosis, for which one could use the critical control parameters of the unperturbed order parameter bifurcation to approximate the optimal control parameters of stochastic resonance, under a suitable scale transform, in designing resonance detecting devices.

## Acknowledgements

This work is supported by the National Natural Science Foundation of China (Grant Nos. 11372233 and 11772241) and by the Hong Kong Research Grants Council under the GRF Grant